\documentclass[bookmarksopen,CJK]{article} 

\usepackage[CJKbookmarks=true,dvipdfm]{hyperref}
\usepackage{amssymb} 
 \usepackage{indentfirst,latexsym, bm}

\begin{document}

  \title{ Divergence form nonlinear nonsmooth elliptic equations with locally arbitrary growth conditions and  nonlinear maximal regularity}
\author{ Qiao-fu Zhang
\\
\vspace{2mm}
\footnotesize{(Academy of Mathematics and Systems Science,}
\\
\footnotesize{
Chinese Academy of Sciences, Beijing 100190, P. R. China)}
     }  

\date{}
\maketitle

\newcommand{\dif}{\,\mathrm{d}}

\newtheorem{dy}{ \textbf{Definition}}[section]
\newtheorem{mt}{ \textbf{Proposition}}[section]
\newtheorem{yl}{ \textbf{Lemma}}[section]
\newtheorem{z}{ \textbf{Remark}}[section]
\newtheorem{dl}{\textbf{Theorem}}[section]
\newtheorem{tl}{\textbf{Corollary}}[section]
\newcommand{\pf}{\textbf{Proof\,\,\,}}
\newcommand{\epf}{\hfill$\square$}

\makeatletter 
\renewcommand\theequation{\thesection%
                   .\arabic{equation}}
\abstract{
             This is a  simplif\mbox{}ication of our prior work on the existence theory for the Rosseland-type equations.
              Inspired  by the   Rosseland equation in the conduction-radiation coupled heat transfer,
               we use the locally arbitrary growth conditions instead of the common global restricted growth conditions.
               In the Lebesgue square integrable space, the solution to the linear elliptic equation   depends     continuously on the coef\mbox{}f\mbox{}icients matrix. This is a simple version of the maximal regularity. There exists a f\mbox{}ixed point for the  linearized map (compact and continuous) in a closed convex set.
}

\vspace{5mm}

\section{Introduction}

  Consider  the following elliptic problem: f\mbox{}ind $u$, $(u-u_b)\in H^1_0(\Omega)$, such that
\begin{equation}
 -\mbox{div}[A(u(x),x)\nabla u]  = 0 ,\quad \mbox{in\,}\,\Omega.
\end{equation}

For the Rosseland equation: $A(z,x)=K(x)+z^3B(x)$, where $K(x)$ and $B(x)$ are symmetric and positive def\mbox{}inite.

   (1)  $K(x)+z^3B(x)$ is  positive def\mbox{}inite only in an interval for $z$.

    (2)  it doesn't satisfy the common growth and smooth conditions and there may be no $C^{2,\gamma}$ estimate (Theorem 15.11 \cite{gl}).

The problem  of the existence theory for the Rosseland equation (also named dif\mbox{}fusion approximation) was proposed by Laitinen \cite{Laitinen} in 2002.
 It may be useful to keep this equation in mind while reading this paper.

It's a little technical to prove the existence by the f\mbox{}ixed point method in $L^\infty(\Omega)$ \cite{amm,zqfthesis}. We will use $L^2(\Omega)$ in this paper.

F\mbox{}irstly, we make the following assumptions.

(A1)  $\Omega\subset \mathbb{R}^n$ is  a bounded Lipschitz domain. 


(A2) $A=(a_{ij})$. $a_{ij} =a_{ji} $.  $T_{min}\leq T_{max}$ are two constants.
\begin{equation}
   \lambda |\xi|^2 \leq a_{ij}(z,x)\xi_i\xi_j\leq \Lambda  |\xi|^2, \quad 0< \lambda  \leq \Lambda ,
  \end{equation}
\begin{equation}
\forall\,\,(z,x,\xi )\in  [T_{min},T_{max} ] \times\Omega  \times \mathbb{R}^n.
\end{equation}

(A3)   \begin{equation}
u_b \in H^{1}( \Omega ).\quad  T_{min}\leq u_b(x)\leq T_{max},\quad \mbox{a.\,e.\,\,in}\,\,\,   \partial\Omega.
\end{equation}

(A4)  $A(z,x )$ is uniformly continuous with respect to $z$ in $\mathfrak{C}$, where
\begin{equation}
 \mathfrak{C} =\{ \varphi\in L^2(\Omega);\,\,T_{min}\leq \varphi(x) \leq T_{max},\,\,\mbox{a.\,e.\,\, in} \,\,\Omega\};
\end{equation}
In other words, if $ z_i,\,z\in \mathfrak{C}$, $\|z_i-z\|_2\to 0$,
  \begin{equation}
  \sup_{1\leq p,q\leq n}\|a_{pq}(z_{i }(x),x)-a_{pq}(z(x),x)\|_2\to 0.
  \end{equation}

\begin{z}\rm
In fact, we had considered a  general case:  parabolic equations with bounded mixed boundary conditions and nonnegative bounded right-hand term $f(z,x)$ in \cite{zqfthesis}.

If $a_{pq}$ is uniformly H\"older continuous    with respect to $z$, (A4) is natural since
  \begin{eqnarray}
  \|a_{pq}(z_{i }(x),x)-a_{pq}(z(x),x)\|_2^2&\leq& \int_{\Omega} C |z_{i }(x)-z(x)  |^{2\alpha}
  \nonumber\\
&\leq& C\|z_i-z\|_2^2\to 0.
  \end{eqnarray}

\end{z}

\section{Linearized map and f\mbox{}ixed point}

\begin{dl} $($Corollary 11.2 \cite{gl}$)$\label{dl:fdp}
Let $\mathfrak{C}$ be a closed convex  set in a Banach space $\mathfrak{B}$ and let $\mathcal{ L}$ be a continuous mapping of $\mathfrak{C}$  into itself such that
the image $\mathcal{ L}\mathfrak{C}$  is precompact. Then $\mathcal{ L}$ has a f\mbox{}ixed point.
\end{dl}
\begin{yl}
The following set
\begin{equation}
 \mathfrak{C} =\{ \varphi\in L^2(\Omega);\,\,T_{min}\leq \varphi(x) \leq T_{max},\,\,a.\,e.\,\, in \,\,\Omega\}
\end{equation}
is a closed convex  set in the Banach space $L^2(\Omega )$.
\end{yl}
\pf Suppose $v_i\in \mathfrak{C}$, $v\in L^2(\Omega )$, $\|v_i-v\|_2\to 0$. If $v \notin \mathfrak{C}$, there exist two constants $\delta_0>0$, $\delta_1>0$, such that the Lebesgue measure of the set
$\Omega_0\equiv\{x\in \Omega;\, v(x)\geq T_{max}+\delta_0\}$ is bigger than $\delta_1>0$. Then
\begin{equation}
 \|v_i-v\|_2^2=\int_\Omega|v_i-v | ^2\geq   \int_{\Omega_0 } |v_i-v | ^2\geq \delta_0^2 \delta_1.
\end{equation}
It's impossible since $\|v_i-v\|_2\to 0$. Similarly, $v \geq T_{min}$ and $\mathfrak{C}$ is closed.
\begin{equation}
   \forall\, \theta\in [0,1],\quad \theta v_1 + (1- \theta) v_2\leq \theta T_{max} + (1- \theta)T_{max}=T_{max}   .
    \end{equation}
     So $\mathfrak{C}$ is convex.
\epf

\begin{dl}
If $(A1)-(A4)$ are satisf\mbox{}ied, then

$(1)$ $\forall\,  z\in \mathfrak{C} $, the following  equation has a unique solution $w \in \mathfrak{C}$, $(w -u_b)\in H^{1}_0( \Omega)$.
\begin{equation}\label{}
\int_\Omega A(z(x),x )\nabla w \cdot\nabla \varphi =0,
\quad \forall\, \varphi\in H^{1}_0( \Omega).
\end{equation}

$(2)$  Def\mbox{}ine a map $\mathcal{ L}:\, \mathfrak{C}\to  \mathfrak{C}$, $\mathcal{ L}z=w$, then $\mathcal{ L}\mathfrak{C}$ is precompact in $L^2(\Omega)$.

$(3)$  $\mathcal{ L}$ is continuous in $L^2(\Omega)$. So $\mathcal{ L}$ has a f\mbox{}ixed point in $\mathfrak{C}$.
\end{dl}
\pf
(1) Let $v=(w -u_b)\in H^{1}_0( \Omega)$, then we have
\begin{equation}
\int_\Omega A(z(x),x )\nabla v \cdot\nabla \varphi =-\int_\Omega A(z(x),x )\nabla u_b \cdot\nabla \varphi,
\quad \forall\, \varphi\in H^{1}_0( \Omega).
\end{equation}

From (A2), if $ z\in \mathfrak{C} $, $A(z(x),x )\in [ \lambda,\Lambda] $. From (A3), $|\nabla u_b|\leq C$.
 Let $ \varphi=v$,
\begin{eqnarray}
\lambda|v|_1^2&\leq& \int_\Omega A \nabla v \cdot\nabla v =-\int_\Omega A \nabla u_b \cdot\nabla v
\nonumber\\
&\leq& \|A \nabla u_b \|_2 \|\nabla v \|_2 =(\int_\Omega A^\top A \nabla u_b \cdot\nabla u_b)^{1/2}\|\nabla v \|_2
\nonumber\\
&\leq&  \Lambda\|u_b\|_{H^{1 }  (\Omega )} |v|_1.
\end{eqnarray}

 Using the well-known Lax-Milgram Lemma, there exists a unique solution $v \in H^{1}_0( \Omega)$ to this equation. Using the Poincar\'e inequality,
\begin{eqnarray}
\|w  \|_{H^{1} ( \Omega)}&\leq& \|w -u_b\|_{H^{1} ( \Omega)}+ \| u_b\|_{H^{1} ( \Omega)}
\nonumber\\
&\leq& C(\Omega)\frac {\Lambda\|u_b\|_{H^{1 }  (\Omega )}}{\lambda }+ \| u_b\|_{H^{1} ( \Omega)}.
\end{eqnarray}

Using the maximum principle (Theorem 8.1 \cite{gl}), $u\in  \mathfrak{C}$. In fact, we can let $ \varphi= (w-T_{max})_+\in H^{1}_0 ( \Omega)$,
\begin{eqnarray}
C(\Omega)\lambda\|(w-T_{max})_+\|_1^2
&\leq& \lambda|(w-T_{max})_+|_1^2
\nonumber\\
&\leq& \int_\Omega A(z(x),x )\nabla (w-T_{max})_+ \cdot\nabla (w-T_{max})_+
\nonumber\\
&=& \int_\Omega A(z(x),x )\nabla w \cdot\nabla (w-T_{max})_+= 0.
\end{eqnarray}
So $(w-T_{max})_+=0$, $w\leq T_{max}$. In the same way, $w\leq T_{min}$.

(2) $\|w  \|_{H^{1} ( \Omega)}\leq C$.  $\mathcal{ L}\mathfrak{C}$ is bounded in $H^{1} ( \Omega)$. From the Rellich Theorem, $\mathcal{ L}\mathfrak{C}$ is precompact in $L^2(\Omega)$.

(3)  Suppose
\begin{equation}
z_i,\,z\in \mathfrak{C},\quad \|z_i-z\|_2\to 0,\quad \mathcal{ L}z_i=w_i , \quad \mathcal{ L}z =w.
\end{equation}

 $ H^{1}_0 ( \Omega)$ is a Hilbert and then a ref\mbox{}lexive space, so there exists a subsequence $\{i_k\}$ and $v_0=(w_0-u_b)\in H^{1}_0 ( \Omega)$ such that
\begin{equation}
(w_{i_k}-u_b) \to (w_0-u_b) ,\quad \mbox{weakly\,\,in\,\,}H^{1}_0 ( \Omega) .
\end{equation}
\begin{equation}
H^{1}_0 ( \Omega)\subset  L^{2}( \Omega) ,\quad (L^{2}( \Omega))'\subset(H^{1}_0 ( \Omega))' .
\end{equation}
\begin{equation}
  (w_{i_k}-u_b) \to (w_0-u_b) ,\quad \mbox{weakly\,\,in\,\,}L^{2}( \Omega) .
\end{equation}

  $\{w_{i_k}-u_b\}$ is bounded in $H^1(\Omega)$, so there exists a subsequence $\{i_m\}\subset \{i_k\}$ and $v_* \in L^{2}( \Omega)$ such that
\begin{equation}
(w_{i_m}-u_b) \to v_* ,\quad \mbox{strongly\,\,in\,\,}L^{2}( \Omega) .
\end{equation}
\begin{equation}
  (w_{i_m}-u_b) \to v_0  ,\quad \mbox{weakly\,\,in\,\,}L^{2}( \Omega) .
\end{equation}
So $v_*=v_0$. Since each subsequence of $  \{\|w_{i_k} -u_b -v_0\|_2\}$ has a sub-subsequence which  converges to 0,  $\|w_{i_k} -u_b -v_0\|_2\to 0 $, $\|w_{i_k} -w_0\|_2\to 0 $.

Since
\begin{equation}
(w_{i_k}-u_b) \to (w_0-u_b) ,\quad \mbox{weakly\,\,in\,\,}H^{1}_0 ( \Omega) .
\end{equation}
  \begin{equation}
\forall\,\overrightarrow{\psi}\in L^2 ( \Omega; \,\mathbb{R}^n),\quad \langle \overrightarrow{\psi}, h\rangle_{H^1_0}\equiv \int_\Omega \nabla h \cdot \overrightarrow{\psi},\quad \forall\, h\in  H^{1}_0 ( \Omega) ,
   \end{equation}
   is a bounded linear functional.
  \begin{equation}
 \int_\Omega \nabla (w_{i_k}-u_b) \cdot \overrightarrow{\psi}\to \int_\Omega  \nabla (w_0-u_b) \cdot \overrightarrow{\psi}.
   \end{equation}
     \begin{equation}
\forall\,\overrightarrow{\psi}\in L^2 ( \Omega; \,\mathbb{R}^n),\quad \int_\Omega \nabla  w_{i_k}  \cdot \overrightarrow{\psi}\to \int_\Omega  \nabla  w_0  \cdot \overrightarrow{\psi}.
   \end{equation}

     From the Riesz representation theorem,  the dual space
  \begin{equation}
  ( L^2 ( \Omega; \,\mathbb{R}^n))'\simeq L^2 ( \Omega; \,\mathbb{R}^n),
   \end{equation}
  \begin{equation}
  \nabla w_{i_k} \to \nabla w_0 ,\quad \mbox{weakly\,\,in\,\,} L^2 ( \Omega; \,\mathbb{R}^n).
   \end{equation}

  From (A4) and $\|z_i-z\|_2\to 0$,
  \begin{equation}
  \sup_{1\leq p,q\leq n}\|a_{pq}(z_{i_k}(x),x)-a_{pq}(z(x),x)\|_2\to 0.
  \end{equation}

  We can conclude that, $\forall\, \phi\in C^{\infty}_0(\Omega)$,
   \begin{eqnarray}
 &&|\int_\Omega [A(z_{i_k} (x),x )\nabla w_{i_k} -A(z  (x),x )\nabla w_0] \cdot\nabla \phi |
\nonumber\\
&\leq& |\int_\Omega [A(z_{i_k} (x),x )\nabla w_{i_k} -A(z  (x),x )\nabla w_{i_k}] \cdot\nabla \phi |
\nonumber\\&&\,+\,|\int_\Omega [A(z  (x),x )\nabla w_{i_k} -A(z  (x),x )\nabla w_0] \cdot\nabla \phi |
\nonumber\\
&=& |\int_\Omega [A(z_{i_k} (x),x )  -A(z  (x),x )] \nabla w_{i_k}\cdot\nabla \phi |
\nonumber\\&&\,+\,|\int_\Omega [\nabla w_{i_k} - \nabla w_0] \cdot A(z  (x),x )^\top\nabla \phi |
\nonumber\\
&\leq& C \sup_{1\leq p,q\leq n}\|a_{pq}(z_{i_k}(x),x)-a_{pq}(z(x),x)\|_2+\epsilon(i_k)
\to 0
 .
\end{eqnarray}
    \begin{equation}
 \int_\Omega A(z_{i_k}(x)  ,x )\nabla w_{i_k}  \cdot\nabla \phi=0,\quad
 \int_\Omega A(z (x)  ,x )\nabla w_0  \cdot\nabla \phi =0
 .
\end{equation}
    \begin{equation}
 \int_\Omega A(z  (x),x )\nabla w_0  \cdot\nabla \varphi =0,
\quad \forall\, \varphi\in H^{1}_0(\Omega).
\end{equation}

Since the solution is unique from the step (1), $w_0=\mathcal{ L}z =w$. So  $\|w_{i_k} -w \|_2\to 0 $.
 Each subsequence of $  \{\|w_{i } -w \|_2 \}$ has a sub-subsequence which  converges to $0$, so $\|w_{i } -w \|_2\to 0 $. We have obtain the continuity of $\mathcal{ L}$.

  From Theorem \ref{dl:fdp}, there exists a f\mbox{}ixed point.
\epf
\begin{z}\rm
For the continuity of $\mathcal{ L}$ in $C^0(\overline{\Omega})$, we can use the well-known De Giorgi-Nash estimate: $\{w_i\}$ is bounded in $C^{0,\alpha} (\overline{\Omega})$ if
$u_b\in C^{0,\alpha_0} ( \partial\Omega )$ and $\Omega$ satif\mbox{}ies a uniform exterior cone condition (Theorem 8.29 \cite{gl}).

Then from the Arzel$\grave{\rm{a}}$-Ascoli Lemma, $\|w_{i_k} -w _0\|_{C^0(\overline{\Omega})}\to 0$. By the same method, $w_0=w$ and $\|w_{i } -w  \|_{C^0(\overline{\Omega})}\to 0$.

From the linear maximal regularity \cite{grell},  a natural conjecture is:    $\mathcal{ L}$ is continuous in $C^{0,\alpha} (\overline{\Omega})$ and $H^1(\Omega)$.
\end{z}

\begin{dy}
\begin{equation}
 \mathfrak{C}_\infty =\{ \varphi\in L^\infty(\Omega);\,\,T_{min}\leq \varphi(x) \leq T_{max},\,\,a.\,e.\,\, in \,\,\Omega\}
\end{equation}
is a closed convex  set in the Banach space $L^\infty(\Omega )$.

We replace $(A4)$ in $L^2$ with $(A4')$ in $L^\infty:$
\begin{equation}
 \|A(z_i(x),x )-A(z(x),x)\|_\infty\to 0, \quad \mbox{if\,\,} \| z_i(x)-z(x) \|_\infty\to 0.
 \end{equation}
\end{dy}

\begin{tl}
If $(A1)-(A3),\,\,(A4')$ are satisf\mbox{}ied,   def\mbox{}ine a map $\mathcal{ L}:\, \mathfrak{C}_\infty\to  \mathfrak{C}_\infty$, $\mathcal{ L}z=w:$ such that $w \in \mathfrak{C}_\infty$, $(w -u_b)\in H^{1}_0( \Omega)$ and
\begin{equation}
\int_\Omega A(z(x),x )\nabla w \cdot\nabla \varphi =0,
\quad \forall\, \varphi\in H^{1}_0( \Omega).
\end{equation}

Then  $\mathcal{ L}$ is continuous in $H^1(\Omega)$.
\end{tl}
\pf  Suppose
\begin{equation}
z_i,\,z\in \mathfrak{C}_\infty,\quad \|z_i-z\|_\infty\to 0,\quad \mathcal{ L}z_i=w_i , \quad \mathcal{ L}z =w.
\end{equation}

(1) $\{w_i\}$ is bounded in $H^1(\Omega)$. For any $\overrightarrow{\psi} \in L^2(\Omega;\, \mathbb{R}^n)$, each subsequence of $ \int_{\Omega }   \nabla (w_i-w)  \cdot \overrightarrow{\psi} $ has a sub-subsequence converges to 0. So
 $ \nabla w_i\to  \nabla w$ weakly in $L^2(\Omega;\, \mathbb{R}^n)$. $ \|A(z_i(x),x )-A(z(x),x)\|_\infty\to 0$,
\begin{equation}
  \nabla u_b \in L^2(\Omega;\, \mathbb{R}^n),\quad  \int_{\Omega }    A(z_i) \nabla w _i \cdot  \nabla u_b
 \to  \int_{\Omega }    A(z) \nabla w  \cdot  \nabla u_b.
\end{equation}
\begin{equation}
  \int_{\Omega }    A(z_i) \nabla w _i \cdot  \nabla u_b = \int_{\Omega }    A(z_i) \nabla w _i \cdot  \nabla w_i,
  \end{equation}
\begin{equation}
    \int_{\Omega }    A(z) \nabla w  \cdot  \nabla u_b=\int_{\Omega }    A(z) \nabla w  \cdot  \nabla w.
\end{equation}
\begin{equation}
 \int_{\Omega }    A(z_i) \nabla w _i \cdot  \nabla w_i\to
 \int_{\Omega }    A(z) \nabla w  \cdot  \nabla w.
\end{equation}
\begin{equation}
   \int_{\Omega }    [A(z_i)-A(z)] \nabla w _i \cdot  \nabla w_i\leq C\|A(z_i(x),x )-A(z(x),x)\|_\infty\to 0.
  \end{equation}
  \begin{equation}
   \int_{\Omega }     A(z) [\nabla w  \cdot  \nabla w- \nabla w _i \cdot  \nabla w_i]\to 0.
  \end{equation}
\begin{eqnarray}
 & & \int_{\Omega }   |\nabla w_{i } - \nabla w  |^2
\nonumber\\
 &\leq & \frac {1}{\lambda} \int_{\Omega }    A(z)[\nabla w_{i } - \nabla w ]\cdot [\nabla w_{i } - \nabla w ]
\nonumber\\
 &= & \frac {1}{\lambda} \int_{\Omega }    A(z)[\nabla w  \cdot  \nabla w+\nabla w_i  \cdot  \nabla w_i-2\nabla w  \cdot  \nabla w_i]
\nonumber\\
&\to& 0.
\end{eqnarray}
$(w_{i } -  w)\in H^1_0( \Omega)$, then use the Poincar\'e inequality, $\|w_i-w\|_{H^1}\to 0$.  \epf

\section{Nonlinear maximal regularity}

For the linear  parabolic/ellptic equations with nonsmooth data, the theory of maximal regularity has been established \cite{gri,gr07,grell,gr}.
In brief, maximal regularity is about the smoothness of the data-to-solution-map \cite{gr}. This smooth dependence has its physical meaning: many physical processes  are stable with respect to the parameters (except the chaos and critical theory).
 For the mathematicians, "the door is open to apply the powerful theorems of dif\mbox{}ferential calculus"(\cite{gr}, e.g.  the Implicit Function Theorem).

In the following, we will discuss the continuous dependence (between the solutions and the data) for the following kind of nonlinear equations: f\mbox{}ind $u$, $(u-u_b)\in H^1_0(\Omega)$, such that
\begin{equation}
\int_\Omega A(u(x),x )\nabla u \cdot\nabla \varphi =0,
\quad \forall\, \varphi\in H^{1}_0( \Omega).
\end{equation}

In this section, we strengthen (A4) with the following equicontinuous conditions.

\textbf{(A4.1)} Each element of $\mathfrak{A}(z,x )$ satisf\mbox{}ies (A2):
\begin{equation}
A(z,x )|_{ \mathfrak{C}\times\Omega }\in [\lambda, \Lambda],\quad \forall\, A(z,x )\in \mathfrak{A}(z,x ).
\end{equation}
\begin{equation}
 \mathfrak{C} =\{ \varphi\in L^2(\Omega);\,\,T_{min}\leq \varphi(x) \leq T_{max},\,\,\mbox{a.\,e.\,\, in} \,\,\Omega\}.
\end{equation}

\textbf{(A4.2)} $\mathfrak{A}(z,x )$ is  equicontinuous with respect to $z$ in $\mathfrak{C}$ (uniformly for $x$). In other words, for any $ z_i,\,z\in \mathfrak{C}$, if $\|z_i-z\|_2\to 0$,
  \begin{equation}
 \sup_{A=(a_{pq})\in \mathfrak{A}(z,x )} \sup_{1\leq p,q\leq n}\|a_{pq}(z_{i }(x),x)-a_{pq}(z(x),x)\|_2\leq \epsilon(i) .
  \end{equation}
$\epsilon(i)$ denotes a higher-order inf\mbox{}initesimal which  depends only on $i$.

\begin{dl}
If $(A1)-(A3),\,\,(A4.1)-(A4.2)$ are satisf\mbox{}ied and the solution to the nonlinear equation is unique, then

$(1)$    $u$ depends continuously $($in $L^2$ or $C^0(\overline{\Omega}))$ on the coef\mbox{}f\mbox{}icients matrix $A(\cdot,x)$ in  $\mathfrak{A}(z,x )$.

$(2)$   $u$ depends continuously $($in $L^2$ or $C^0(\overline{\Omega}))$ on the boundary value $u_b$.

$(3)$  $A(z,x)\nabla z$ is a continuous $($with respect to $z$ in $\mathfrak{C}\cap H^1(\Omega))$ functional in $(C^{0,1}(\overline{\Omega}))'$.
\end{dl}
\pf
(1)  Suppose $A_i(\cdot,x)\in \mathfrak{A}(\cdot ,x )$, then there exists a $u_i$, $(u_i-u_b)\in H^1_0(\Omega)$, such that
\begin{equation}
\int_\Omega A_i(u_i(x),x )\nabla u_i \cdot\nabla \varphi =0,
\quad \forall\, \varphi\in H^{1}_0( \Omega).
\end{equation}

We will prove that $\|u-u_i\|_{2 }\to 0$ (or $\|u-u_i\|_{C^0(\overline{\Omega}) }\to 0$)  if $\|A(z(x),x)- A_i(z(x),x)\|_{ 2}\to 0$ for any $z(x)\in \mathfrak{C} $.
 $\{u_i\}$ is bounded in $H^1(\Omega)$ (or in $C^{0,\alpha}(\overline{\Omega})$).
  So 　$u_{i_k}\to u_0$, strongly in $L^2(\Omega)$ (or in $C^{0}(\overline{\Omega})$);  $ \nabla u_{i_k}\to  \nabla u$, weakly in $L^2(\Omega;\, \mathbb{R}^n)$.
 \begin{eqnarray}
 &&\| A_{i_k}(u_{i_k}(x),x )  -A(u_0(x),x )\|_2
\nonumber\\
 &\leq &
\| A_{i_k}(u_{i_k}(x),x ) -A_{i_k}(u_0(x),x )\|_2+\| A_{i_k}(u_0(x),x ) -A(u_0(x),x )\|_2
\nonumber\\
 &\to &  0.
\end{eqnarray}

 $\forall\,\phi\in C^{0,\infty}(\Omega)$,
\begin{equation}
0=\int_\Omega A_{i_k}(u_{i_k}(x),x )\nabla u_{i_k} \cdot\nabla \phi \to  \int_\Omega A (u_0(x),x )\nabla u_0 \cdot\nabla \phi.
\end{equation}

Because of the density,
 \begin{equation}
\int_\Omega A (u_0(x),x )\nabla u_0 \cdot\nabla  \varphi =0,
\quad \forall\, \varphi\in H^{1}_0( \Omega).
\end{equation}

Since the solution is unique, so $w_i\to w_0$, strongly in $L^2(\Omega)$ (or in $C^{0}(\overline{\Omega})$).

(2) Let $u_{bi},\,u_{b0}\in H^1(\Omega)$, $\|u_{bi}-u_{b0}\|_{ H^1(\Omega)}\to 0$. $w_i=u_i-u_{bi}\in H^1_0(\Omega)$ such that $ \forall\, \varphi\in H^{1}_0( \Omega)$,
 \begin{equation}
\int_\Omega A (w_i+u_{bi} )\nabla w_i \cdot\nabla  \varphi =-\int_\Omega A (w_i+ u_{bi} )\nabla u_{bi} \cdot\nabla  \varphi.
\end{equation}

 $\{w_i\}$ is bounded in $H^1(\Omega)$ (or in $C^{0,\alpha}(\overline{\Omega})$).
  So 　$w_{i_k}\to w_0$, strongly in $L^2(\Omega)$ (or in $C^{0}(\overline{\Omega})$);  $ \nabla w_{i_k}\to  \nabla w_0$, weakly in $L^2(\Omega;\, \mathbb{R}^n)$.
By the same method in (1), $\forall\,\phi\in C^{0,\infty}(\Omega)$,
 \begin{equation}
\int_\Omega A (w_0+u_{b0} )\nabla w_0 \cdot\nabla  \phi =-\int_\Omega A (w_0+ u_{b0} )\nabla u_{b0} \cdot\nabla  \phi.
\end{equation}

  So 　$w_{i }\to w_0$, strongly in $L^2(\Omega)$.

  (3) For any $\eta\in C^{0,1}(\overline{\Omega})$,
   \begin{eqnarray}
|\langle A(z,x)\nabla z, \eta\rangle|
&=&|\int_\Omega A(z,x)\nabla z \cdot\nabla \eta|
\nonumber\\
 &\leq &
  \|\eta\|_{ C^{0,1}(\overline{\Omega})}\int_\Omega |A(z,x)\nabla z|
\nonumber\\
 &\leq &
  C\|\eta\|_{ C^{0,1}(\overline{\Omega})}.
\end{eqnarray}

For any $z$ in $\mathfrak{C}\cap H^1(\Omega)$,  $\langle A(z,x)\nabla z, \eta\rangle$ is a linear continuous   functional in $(C^{0,1}(\overline{\Omega}))'$.

Suppose $z_i\to z$ in $H^1(\Omega)$,
   \begin{eqnarray}
&&|\langle A(z_i,x)\nabla z_i- A(z,x)\nabla z , \eta\rangle|
\nonumber\\
 &\leq&  \|\eta\|_{ C^{0,1}(\overline{\Omega})} |\int_\Omega |A(z_i,x)\nabla z_i- A(z,x)\nabla z|
 \to 0.
\end{eqnarray}

$\| A(z_i,x)\nabla z_i- A(z,x)\nabla z\|_{ (C^{0,1}(\overline{\Omega}))'}\to 0$. \epf

\begin{z}\rm
We can consider the continuous dependence in $H^1(\Omega)$.
\end{z}

\section{Existence for Garlerkin method }

 Let $h\in (0,1)$ be the step size,  $\{\phi_{i,h } \}  $ is a kind of f\mbox{}inite element basis in  $H^1_0(\Omega)$.

  $\forall\,  z_h=(\sum_i z_{i,h}\phi_{i,h }+u_b)\in \mathfrak{C} $, the following  equation has a unique solution $w_h=\sum_i w_{i,h}\phi_{i,h }$, $(w_h+u_b) \in \mathfrak{C}\cap H^{1}_0( \Omega)$.
\begin{eqnarray} 
&&\int_\Omega A(\sum_i z_{i,h}\phi_{i,h }+u_b)\nabla (\sum_i w_{i,h}\phi_{i,h }) \cdot\nabla \phi_{j,h }
 \nonumber\\
 &=& 
 -
\int_\Omega A(\sum_i z_{i,h}\phi_{i,h }+u_b)\nabla u_b\cdot\nabla \phi_{j,h } 
,
\quad \forall\, \phi_{j,h }.
\end{eqnarray} 

Def\mbox{}ine $\mathcal{L}z_h =(w_h+u_b) $, then 

(1) $\mathcal{L} \mathfrak{C}\subset \mathfrak{C} $.

(2) $\mathcal{L}\mathfrak{C}$ is compact.

(3)  $\mathcal{L} $ is continuous.

(4) There exists a f\mbox{}ixed point $u_{h}=\sum_i w_{i,h}\phi_{i,h }$ and $u_h\to u$.

\section{Acknowledge }

This work is supported by the National Nature Science Foundation of China (No. 90916027).
This is a part  of my PhD thesis \cite{zqfthesis} in AMSS,  Chinese Academy of Sciences, and a simplif\mbox{}ication of our prior paper \cite{amm}.  So I will  thank my advisor Professor Jun-zhi Cui (he is also a member of the Chinese Academy of Engineering) and the referees for their careful reading and helpful comments. My E-mail is:  zhangqf@lsec.cc.ac.cn.

%
%

\end{document}